\documentclass[a4paper,11pt]{amsart}
\usepackage{amssymb}
\usepackage{amsthm, amsmath, latexsym, amsfonts}
\usepackage[all]{xy}
\newtheorem{Lemma}{Lemma}

\newtheorem{Theorem}{Theorem}

\newtheorem{Question}{Question}
\theoremstyle{definition}
\newtheorem{Remark}[Theorem]{Remark}

\newcommand{\C}{\mathbb{C}}

\newcommand{\mgn}{\mathcal{M}_{g,n}}

\newcommand{\mgnbar}{\overline{\mathcal{M}}_{g,n}}

\title[Holomorphic forms on moduli spaces of stable curves]{Holomorphic differential forms \\ on moduli spaces of stable curves}
\author{Claudio Fontanari}

\email{claudio.fontanari@unitn.it}\curraddr{
Universit\`a degli Studi di Trento \\
Dipartimento di Matematica \\
Via Sommarive 14 \\
38123 Povo (Trento) \\ Italy.}

\begin{document}

\begin{abstract}
We prove that the space of holomorphic $p$-forms on the moduli space 
$\overline{\mathcal{M}}_{g,n}$ of stable curves of genus $g$ with $n$ 
marked points vanishes for $p=14, 16, 18$ unconditionally and also for $p=20$ 
under a natural assumption in the case $g=3$. This result is consistent with the 
Langlands program and it is obtained by applying the Arbarello-Cornalba 
inductive approach to the cohomology of moduli spaces. 
\end{abstract}

\maketitle

\section{Introduction}

The moduli space $\mgnbar$ parameterizing stable curves of genus $g$ with $n$ 
marked points is a projective compactification with a beautiful geometric structure:
all its boundary components are (products of) moduli spaces of the same 
kind but with smaller invariants. This remarkable property was employed by Enrico 
Arbarello and Maurizio Cornalba to perform an elegant inductive computation of the first few 
rational cohomology groups of $\mgnbar$. In particular, in \cite{AC} they proved 
that $H^1(\mgnbar) = H^3(\mgnbar) = H^5(\mgnbar) = 0$ and established an 
inductive approach to reduce the vanishing of odd cohomology (so long as it 
vanishes, since it is well known that $H^{11,0}(\overline{\mathcal{M}}_{1,11}) \ne 0$)
to a finite number of explicit verifications in low genus. 

A few years later, in \cite{BF} Gilberto Bini and the author pointed out that the same 
inductive procedure implies also the vanishing of the spaces of holomorphic $p$-forms 
$H^{p,0}(\mgnbar)$ for $0 < p < 11$. 
More recently, a renewed interest in the Arbarello-Cornalba method is witnessed 
by the papers \cite{BFP} by Jonas Bergstr\"om, Carel Faber, and Sam Payne, where 
they compute that $H^7(\mgnbar) = H^9(\mgnbar) = 0$, and \cite{CLP} by Samir Canning, 
Hannah Larson, and Sam Payne, where they prove inductively that the cohomology group 
$H^k(\mgnbar)$ is pure Hodge-Tate (hence, in particular, $H^{k,0}(\mgnbar)=0$) 
for any even $k \le 12$. This is consistent with the Langlands program,
predicting that $H^k(\mgnbar)$ should be pure Hodge-Tate for all even $k \le 20$.  

Here we move a small step forward along the same path by obtaining the following result: 

\begin{Theorem}\label{main}
We have 
$$ 
H^{14,0}(\mgnbar) =  H^{16,0}(\mgnbar) =  H^{18,0}(\mgnbar) = 0
$$
for every $g$ and $n$ with $2g-2 + n > 0$. 

\noindent
Furthermore, 
if $H^{20,0}(\overline{\mathcal{M}}_{3,15}) = H^{20,0}(\overline{\mathcal{M}}_{3,16}) = 0$
then $H^{20,0}(\mgnbar) = 0$ for every $g$ and $n$ with $2g-2 + n > 0$. 
\end{Theorem}

Once again, the crucial ingredient is a minor variant of the Arbarello-Cornalba inductive approach 
(see Lemma \ref{hodge}). Of course, the statement of Theorem \ref{main} arises the following natural question: 

\begin{Question}
Is $H^{20,0}(\overline{\mathcal{M}}_{3,15})  =  H^{20,0}(\overline{\mathcal{M}}_{3,16}) = 0$?
\end{Question}

We work over the complex field $\C$.

\medskip
 
{\bf Acknowledgements:} The author is member of GNSAGA of the Istituto Nazionale di Alta Matematica ``F. Severi". This research project was partially supported by PRIN 2017 ``Moduli Theory and Birational Classification''. 

\section{The proofs}

\begin{Lemma}\label{hodge}
Let $0 < p \le 21$ and assume $h^{p,0}(\overline{\mathcal{M}}_{g',n'})=0$ for every $g', n'$ 
such that $p \ge 2g'-2+n' > 0$. Then $h^{p,0}(\mgnbar)=0$ for every $g$ and $n$ 
with $2g-2 + n > 0$. 
\end{Lemma}

\proof 
By double induction on $g$ and $n$. Let $d(g,n) = 2g-2+n > 0$. 

If $d(g,n)=1$ we have either $g=0$ and $n=3$,
or $g=1$ and $n=1$, and in both cases the claim is obvious. 

Let now $d(g,n)>1$. If $p \ge d(g,n)$ then the claim holds by assumption, 
hence let $p < d(g,n)$.
In the long exact sequence of cohomology with compact supports:
$$
\ldots \to H^k_c(\mgn) \to H^k(\mgnbar) \to H^k(\partial \mgnbar) \to \ldots 
$$ 
we have $H^k_c(\mgn) =0$ for $k < d(g,n)$ by \cite{H}.
Since the morphism 
$$
H^k(\mgnbar) \rightarrow H^k(\partial \mgnbar)
$$ 
is compatible with the Hodge structures (see \cite{AC},
p. 102), for $p < d(g,n)$ there is an injection
\begin{equation}
\label{injection}
H^{p,0}(\mgnbar) \hookrightarrow H^{p,0}(\partial \mgnbar).
\end{equation}

Next we use the fact that each irreducible component of the boundary $\partial \mgnbar$
is the image of a map from ${\overline {\mathcal M}}_{g-1, n+2}$ or 
${\overline{\mathcal M}}_{h,m+1} \times {\overline{\mathcal M}}_{g-h, n-m+1}$, 
where $0 \leq h \leq g$ and both $2h-2+m+1$ and $2(g-h)-2+n-m+1$ are
positive. By the analogue of Lemma (2.6) in \cite{AC} and the 
Hodge-K\"{u}nneth formula, the map
\begin{eqnarray*}
 H^{p,0}(\mgnbar) \rightarrow H^{p,0}({\overline {\mathcal M}}_{g-1,n+2}) 
 \oplus \bigoplus_{h, m} H^{p,0}({\overline{\mathcal M}}_{h, m+1} \times
 {\overline{\mathcal M}}_{g-h, n-m+1}) \\
= H^{p,0}({\overline {\mathcal M}}_{g-1,n+2}) 
\oplus \bigoplus_{h, m} (H^{0,0}({\overline{\mathcal M}}_{h, m+1}) \otimes 
H^{p,0}({\overline{\mathcal M}}_{g-h, n-m+1}) \oplus \\
\bigoplus_{q \ge 1} H^{q,0}({\overline{\mathcal M}}_{h, m+1}) \otimes 
H^{p-q,0}({\overline{\mathcal M}}_{g-h, n-m+1}))
\end{eqnarray*} 
is injective whenever the map \eqref{injection} is. 
The right hand side involves the terms
$H^{p,0}({\overline {\mathcal M}}_{g-1,n+2})$
and $H^{p,0}({\overline{\mathcal M}}_{g-h, n-m+1})$ 
with either $h \ge 1$ or $h=0$ and $m \ge 2$,
hence vanishing by induction, 
and products of two terms which have $1 \le q \le 10$, since $p \le 21$. 
Therefore by \cite{BF}, Theorem 1, stating that $H^{q,0}(\mgnbar) =0$ for $0 < q < 11$, 
we obtain $H^{p,0}(\mgnbar)=0$.

\qed

\begin{Remark}
The assumption of Lemma \ref{hodge} is not satisfied for every $11 \le p \le 21$: 
in particular, as it is well known $H^{11,0}(\overline{\mathcal{M}}_{1,11}) \ne 0$ (see for 
instance \cite{FP}, Section 2.3) and  also $H^{17,0}(\overline{\mathcal{M}}_{2,14}) \ne 0$
(see \cite{FP}, Section 3.5).
\end{Remark}

\medskip

\noindent \textit{Proof of Theorem \ref{main}.} In order to apply Lemma \ref{hodge}
we have to fix an even integer $p$ with $14 \le p \le 20$ and check that 
$H^{p,0}(\overline{\mathcal{M}}_{g',n'})=0$ for every $g', n'$ such that 
$p \ge 2g'-2+n' > 0$.

If $g'=0$ then all cohomology is tautological (hence algebraic) by \cite{K}. 

If $g'=1$ then all even cohomology is tautological by \cite{P1}.

If $g'=2$ then all even cohomology is tautological for $n' < 20$ by \cite{P2}.

If $g'=3$ then $\overline{\mathcal{M}}_{g',n'}$ is unirational 
(hence $H^{p,0}(\overline{\mathcal{M}}_{g',n'})=0$ for every $p > 0$)
for $n' \le 14$ by \cite{L}, Theorem 7.1 (notice that this range completely 
covers the case $p \le 18$, while for $p=20$ we need the additional assumption 
in the statement). 

The same Theorem 7.1 in \cite{L} yields the unirationality of $\overline{\mathcal{M}}_{g',n'}$
also for $g'=4$ and $n' \le 15$, $g'=5$ and $n' \le 12$, $g'=6$ and $n' \le 15$, 
$g'=7$ and $n' \le 11$, $g'=9$ and $n' \le 8$, $g'=11$ and $n' \le 10$. 

Finally, by \cite{BCF}, Theorem B., $\overline{\mathcal{M}}_{g',n'}$ is unirational for
$g'=8$ and $n' \le 11$ and $g'=10$ and $n' \le 3$, thus covering the last missing cases. 

\qed


\begin{thebibliography}{99}

\bibitem{AC} E. Arbarello and M. Cornalba: Calculating
cohomology groups of moduli spaces of curves via algebraic
geometry. Inst. Hautes \'Etudes Sci. Publ. Math. 88 (1998), 97--127.

\bibitem{BCF} E. Ballico, G. Casnati, C. Fontanari: On the birational geometry of moduli spaces of pointed curves. Forum Math. 21 (2009), 935--950. 

\bibitem{BFP} J. Bergstr\"om, C. Faber, S. Payne: Polynomial point counts 
and odd cohomology vanishing on moduli spaces of stable curves. 
arXiv:2206.07759 (2022). 

\bibitem{BF} G. Bini and C. Fontanari: Moduli of curves and spin structures via algebraic geometry. Trans. Amer. Math. Soc. 358 (2006), 3207--3217. 

\bibitem{CLP} S. Canning, H. Larson, S. Payne: The eleventh cohomology group of $\mgnbar$. 
arXiv:2209.03113 (2022).

\bibitem{FP} C. Faber and R. Pandharipande: Tautological and non-tautological cohomology of the moduli space of curves. Handbook of moduli. Vol. I, 293--330, Adv. Lect. Math. (ALM), 24, Int. Press, Somerville, MA, 2013. 

\bibitem{H} J. Harer: The virtual cohomological dimension of the
mapping class group of an orientable surface. Invent. Math. 84,
157--176 (1986).

\bibitem{K} S. Keel: Intersection theory of moduli space of stable
$n$-pointed curves of genus zero. Trans. Amer. Math. Soc. 330
(1992), 545--574.

\bibitem{L} A. Logan: The Kodaira dimension of moduli spaces of
curves with marked points. Amer. J. Math. 125, 105--137 (2003). 

\bibitem{P1} D. Petersen: The structure of the tautological ring in genus one. 
Duke Math. J. 163 (2014), 777--793.

\bibitem{P2} D. Petersen: Tautological rings of spaces of pointed genus two 
curves of compact type. Compos. Math. 152 (2016), 1398--1420.

\end{thebibliography}
\end{document}